\documentclass{amsart}
\usepackage{latexsym}
\usepackage{amsmath,amsfonts,epsfig, graphics, graphicx, amsthm, amssymb,mathrsfs}
\input{xy}
\xyoption{all}

\newtheorem{theorem}{Theorem}[section]

\begin{document}

\title[Trilinear Kloosterman fractions II: subdyadic intervals]{Trilinear Kloosterman fractions II: subdyadic intervals and nearly balanced convolutions}
\author[T. Wright]{Thomas Wright}
\begin{abstract}  
This paper broadens the range on which Fouvry and Radziwiłł's results on nearly balanced convolutions apply.  In particular, let $\alpha_m$ and $\beta_n$ be sequences supported on $m\sim M$ and $n\sim N$ where $\beta_n$ is equidistributed for small moduli, and let $Q=X^{\frac 12+\varepsilon}$.  We find that
\begin{gather*}\sum_{q\sim Q}\left|\mathop{\sum\sum}_{\substack{n\sim N,m\sim M \\ mn\equiv a\pmod q}}\alpha_m\beta_n-\frac{1}{\phi(q)}\mathop{\sum\sum}_{\substack{n\sim N,m\sim M \\ (mn,q)=1}}\alpha_m\beta_n\right|\ll \frac{X}{\log^A X}
\end{gather*}
if $N=X^{\frac 12+\delta}$ and $M=X^{\frac 12-\delta}$ with $0<\delta<\frac 1{68}$, which improves Fouvry and Radziwiłł's $0<\delta<\frac 1{112}$.  To prove this, we sharpen Bettin and Chandee's famous result on trilinear forms with Kloosterman fractions in the case where some of the sums are over subdyadic intervals.
\end{abstract}
\maketitle

\section{Introduction}

As is so often the case with papers of this type, we begin with the following question: given an arithmetic function $f$ and an integer $a\neq 0$, what is the largest $\delta>0$ for which $Q\leq X^{\frac 12+\delta}$ implies that
$$\sum_{\substack{q\sim Q \\ (a,q)=1}}\left|\sum_{\substack{n\sim X \\ n\equiv a\pmod q}}f(n)-\frac{1}{\phi(q)}\sum_{\substack{n\sim X \\ (n,q)=1}}f(n)\right|\ll \frac{X}{\log^A X}$$
for some (or possibly for every) $A>0$?  

The motivation for the current paper is Fouvry and Radziwiłł's result on nearly balanced convolutions \cite{FR2}.  In that paper, the authors find the following:
\begin{theorem} \label{FR2}
  Let $k \geq 1$ be an integer and $M,N \geq 1$ be given. Set $X = M N$. Let $\alpha_m$ and $\beta_n$ be two sequences of real numbers supported respectively on $[M, 2M]$ and $[N, 2N]$. Suppose that $\boldsymbol \beta = (\beta_n)$ is Siegel--Walfisz and suppose that $|\alpha_{m}| \leq \tau_k(m)$ and $|\beta_n| \leq \tau_k(n)$ for all integers $m,n \geq 1$. Then, for every $\varepsilon > 0$ and every $A > 0$,
  \begin{equation} \label{eq:maineq}
  \sum_{\substack{Q \leq q \leq 2Q \\ (q,a) = 1}} \Big | \sum_{\substack{m n \equiv a \pmod{q}}} \alpha_m \beta_n - \frac{1}{\varphi(q)} \sum_{(m n, q) = 1} \alpha_m \beta_n \Big | \ll_{A} X (\log X)^{-A} 
  \end{equation}
  uniformly in $N^{56/23} X^{-17/23 + \varepsilon} \leq Q \leq N X^{-\varepsilon}$ and $1 \leq |a| \leq X$.
\end{theorem}
This result, along with those authors' previous result on unbalanced convolutions \cite{FR}, establish greater than 1/2 level of distribution over specific ranges of $N$.  Here, in the case where $Q=X^{\frac 12+\varepsilon}$, Theorem \ref{FR2} means that we can take $N\leq X^{\frac{57}{112}-\varepsilon}$ or $Q^{\frac{57}{56}-\varepsilon}$.

In order to prove Theorem \ref{FR2}, the authors depend on the well-known bound of Bettin and Chandee for trilinear Kloosterman fractions. \cite{BC}.  Define
$$\mathcal B(M,N,A)=\mathop{\sum\sum\sum}_{\substack{a\sim A, m\sim M,n\sim N \\ (m,n)=1}}\alpha_m\beta_n\nu_a e\left(\vartheta \frac{a\overline{m}}{n}\right),$$
where $\vartheta$ is a nonzero integer.

In \cite{BC}, the authors then prove the following.
\begin{theorem}\label{BCThm}
\begin{align*}
\mathcal B(M,N,A)\ll &\|\alpha\|\|\beta\|\|\nu\|\left(1+\frac{|\vartheta|A}{MN}\right)^\frac 12\\
&\times\left((AMN)^{\frac 7{20}+\varepsilon}(M+N)^\frac 14+(AMN)^{\frac 38+\varepsilon}(AN+AM)^\frac 18\right).
\end{align*}
\end{theorem}
This bound is applied at a key juncture in the Fouvry-Radziwiłł paper \cite[(41)]{FR2} to bound a trilinear sum.

\section{New ideas: subdyadic intervals}

In this paper, we aim to broaden (slightly) the allowable range of $N$.  In our previous paper in this series \cite{WrKl}, we slightly enlarged the range of $N$ in Fouvry and Radziwiłł's unbalanced convolutions paper \cite{FR} by improving Theorem \ref{BCThm} for a specific instance that would be helpful to the unbalanced case, namely where the Kloosterman sum has a small fixed term in the denominator.  Here, we embark on a similar approach, improving Theorem \ref{FR2} by improving Theorem \ref{BCThm} in the case where the sums of $m$ and $n$ run over subdyadic intervals.




For intervals $\mathcal A$, $\mathcal M$, and $\mathcal N$, define
$$\mathcal B(\mathcal M,\mathcal N,\mathcal A)=\mathop{\sum\sum\sum}_{\substack{a\in \mathcal A, m\in \mathcal M,n\in \mathcal N \\ (m,n)=1}}\alpha_m\beta_n\gamma_a e\left(\vartheta \frac{a\bar m}{n}\right),$$
where $\vartheta$ is a nonzero integer.  


\begin{theorem}\label{BCsubdy}
Let $\mathcal A\subset [A,2A]$, $\mathcal M\subset [M,2M]$, and $\mathcal N\subset [N,2N]$ be such that for each set, the elements of that set either comprise an interval or are consecutive elements of a congruence class for some modulus.  For some $\eta>0$, assume $|\mathcal M|\ll MX^{-\eta}$ and $|\mathcal N|\ll NX^{-\eta}$.  Then
\begin{align*}
 \mathcal{B}(M,N,A)&\ll \|\alpha\|\|\beta\|\|\nu\|M^\varepsilon \left(1+ \frac{|\vartheta |A}{NM}\right)^\frac12 \\
&\times\hspace{0em}\left(A^\frac 12\left(M^\frac 12N^{\frac 38}+M^\frac 38N^{\frac 12}\right)+A^\frac7{20} \left( M^\frac35N^{\frac7{20}}X^{-\frac{2\eta}{5}}+M^{\frac7{20}}N^{\frac35}X^{-\frac{2\eta}{5}}\right)\right).
\end{align*}
\end{theorem}
In the case where $\alpha$ and $\beta$ are 1-bounded, we also have
$$\|\alpha\|\|\beta\|\ll \sqrt{MN}X^{-\eta}.$$
Since it is generally the case in applications that 
$$A^\frac7{20} \left( M^\frac35N^{\frac7{20}}X^{-\frac{2\eta}{5}}+M^{\frac7{20}}N^{\frac35}X^{-\frac{2\eta}{5}}\right)\gg A^\frac 12\left(M^\frac 12N^{\frac 38}+M^\frac 38N^{\frac 12}\right),$$
our bound for a subdyadic trilinear Kloosterman sum then gives $X^{-\frac{7\eta}{5}}$ times the bound for the full dyadic sum.

The reason that this theorem is helpful for Fouvry-Radziwiłł result \cite{FR2} is as follows.  In \cite[(34)]{FR2}, if we let $\nu_1=\frac{n}{(n,m)}$ and $\nu_2=\frac{n}{(n,m)}$, the authors define the variable $r$ such that
$$\nu_1-\nu_2=qr,$$
and they then replace $q$ in various expressions with $(\nu_1-\nu_2)/r$.  The authors later consider a sum running over $\nu_1$, $\nu_2$, and $r$ that they allow to be unencumbered by all but the most trivial bounds and the congruence requirement that $\nu_1\equiv \nu_2\pmod r$.  When $N$ is close to $Q$ (as is the case in the original paper), the loose bounds on $\nu_1$, $\nu_2$, and $r$ do not make much of a difference.  However, when $N$ is larger than $Q$, $\nu_i$ can be quite a bit larger than $2rQ$, and hence this allows for $q=(\nu_1-\nu_2)/r$ that are much larger than $2Q$. 

To remedy this, we restrict $\nu_1$ and $\nu_2$ to intervals of size $20rQ$, where $\nu_1\equiv \nu_2\equiv j\pmod r$ for $j$ running over the reduced residue classes modulo $r$.  For a given $j$, this means that the $\nu_i$ are summed over sets of size $\ll QX^\varepsilon$.  Letting $X^{-\eta}=N/Q$, we then apply our result to garner a savings from Theorem \ref{FR2}.

Hence, Theorem \ref{BCsubdy} will allow us to prove the following.


\begin{theorem} \label{FRch}
  Let $\alpha_m$, $\beta_n$, $M$, $N$, and $X$ be as in Theorem \ref{FR2}.  Then, for every $\varepsilon > 0$ and every $A > 0$,
  \begin{equation} \label{eq:maineq2}
  \sum_{\substack{Q \leq q \leq 2Q \\ (q,a) = 1}} \Big | \sum_{\substack{m n \equiv a \pmod{q}}} \alpha_m \beta_n - \frac{1}{\varphi(q)} \sum_{(m n, q) = 1} \alpha_m \beta_n \Big | \ll_{A} X (\log X)^{-A} 
  \end{equation}
  uniformly in $N^{34} X^{-17+ \varepsilon} \leq Q \leq N X^{-\varepsilon}$ and $1 \leq |a| \leq X$.

\end{theorem}
In particular, if $Q=X^{\frac 12+\varepsilon}$, then $N$ can be as large as $Q^{\frac{35}{34}-\varepsilon}$ or $X^{\frac{35}{68}-\varepsilon}$.


\section{Proof of Theorem \ref{BCsubdy}: subdyadic intervals and trilinear Kloosterman fractions}

We begin first with some tweaks of the Bettin-Chandee result when the intervals for $m$ and $n$ are subdyadic.  Our change here will be very simple: in instances where there is a trivial sum over $m$ or $n$, we replace the resulting $M$ or $N$ with $MX^{-\eta}$ or $NX^{-\eta}$.

First, in the case of the diagonal terms, we replace any of the trivial bounds of $M$ with $MX^{-\eta}$, and hence \cite[(3.2)]{BC} is now
$$\mathscr{D}_b(M,N,A,L;\beta,\nu) \ll \|\beta\|^2\|\nu\|^2 L\Big(A(bLN)^{\frac 12}+\frac{AM}{bNX^{\eta}}+\frac{M}{X^{\eta}}\Big)M^{\varepsilon}.$$
Next, in the computation of $\mathscr{V}^*_{b,\xi}$, there is a sum of $n_1\in \mathcal N$ that is trivially bounded, and hence we adjust \cite[(4.14)]{BC} appropriately.

$$\mathscr{V}^*_{b,\xi}\ll \frac{\|\nu\|^2b^{\frac12}L^2N^{\frac32}X^{-\eta}M^\varepsilon DA^{2}}{\mathfrak q_1\mathfrak p_2(\mathfrak p_1+\mathfrak q_1)(\mathfrak p_2+\mathfrak q_2)}$$

The computation of $\mathcal V^{\Delta=0}_{b,\xi}$ has two trivially bounded sums over $\mathcal N$, and hence the bound in \cite[(4.22)]{BC} can be replaced with
$$\mathcal V^{\Delta= 0}_{b,\xi}\ll \frac{ \|\nu\|^2AL^2N^{2}X^{-2\eta}M^\varepsilon}{\mathfrak p_1\mathfrak p_2\mathfrak q_1\mathfrak q_2(\mathfrak p_1+\mathfrak q_1)(\mathfrak p_2+\mathfrak q_2)}.$$

The sum $\mathcal V^{\Delta\neq0}_{b,\xi}$ in \cite[(4.29)]{BC} also has a sum over $\mathcal N$ that is trivially bounded, and hence the resulting bound in (4.29) can be replaced with
$$\mathscr V^{\Delta\neq 0}_{b,\xi}\ll \frac{\|\nu\|^2A^2D^2L^5N^{\frac32}X^{-\eta}M^\varepsilon}{(\mathfrak p_1+\mathfrak q_1)(\mathfrak p_2+\mathfrak q_2)\mathfrak p_1\mathfrak q_1^2\mathfrak p_2^3\mathfrak q_2}\left(1+ \frac{|\vartheta|AD}{bLN^2}\right).$$
Putting these together, \cite[(4.20)]{BC} has
$$\mathscr V_{b,\xi}=\mathscr V^{\Delta= 0}_{b,\xi}+\mathscr V^{\Delta\neq 0}_{b,\xi},$$
which gives
$$\mathscr V_{b,\xi}\ll \|\nu\|^2M^\varepsilon\left(1+ \frac{|\vartheta|AD}{bLN^2}\right)\frac{A^2D^2L^{5}N^{\frac32}X^{-\eta}+AL^2N^{2}X^{-2\eta} }{\mathfrak q_1\mathfrak p_2(\mathfrak p_1+\mathfrak q_1)(\mathfrak p_2+\mathfrak q_2)}.$$
We then have
$$\mathscr{U}_{b,\xi}=\mathscr{V}_{b,\xi}+\mathscr{V}^*_{b,\xi},$$
and thus 
$$\mathscr U_{b,\xi} \ll \|\nu\|^2A^2L^{2}N^{\frac32}X^{-\eta}M^\varepsilon
\frac{  D  b^{\frac12}+L^3D^2+N^{\frac12}X^{-\eta}A^{-1}}{\mathfrak q_1\mathfrak p_2(\mathfrak p_1+\mathfrak q_1)(\mathfrak p_2+\mathfrak q_2)}\left(1+ \frac{|\vartheta|AD}{bLN^2}\right).$$
Since 
$$\mathscr{S}_{b,\xi} ^2 \ll M^\varepsilon\|\beta\|^4\|\nu\|^2 \mathop{\sum\sum\sum\sum}\limits_{\substack{\mathfrak p_1,\mathfrak p_2,\mathfrak q_1,\mathfrak q_2\in\mathcal L\cup\{1\},\\ \mathfrak p_1\neq\mathfrak q_1 \Rightarrow 1\in\{\mathfrak p_1,\mathfrak q_1\},\\ \mathfrak p_2\neq\mathfrak q_2 \Rightarrow 1\in\{\mathfrak p_2,\mathfrak q_2\}}}b\mathfrak q_1\mathfrak p_2\left(\mathscr U_{b,\xi} +\mathscr U^*_{b,\xi} \right),$$
and the bound for $\mathscr U^*_{b,\xi}$ is the same as the bound for $\mathscr U_{b,\xi}$, this changes \cite[(4.31)]{BC} to
\begin{gather}\label{ees}
\mathscr{S}_{b,\xi} \hspace{0em}\ll \|\beta\|^2\|\nu\|^2 \Big(b+ \frac{|\vartheta|A}{NM}\Big)^\frac12ALN^{\frac34}X^{-\frac \eta 2}M^\varepsilon\left(\frac{b^{\frac14} N^\frac12 L^\frac12}{M^\frac12}+\frac{L^\frac52N}M+\frac{N^{\frac14}}{ A^\frac12X^{\frac \eta 2}}\right).
\end{gather}
The estimate in \cite[(4.32)]{BC} for $\mathscr{S}^*_{b,\xi}$ is the same.  Since $\mathscr{O}_{b}$ is defined in \cite[(4.1)]{BC} and \cite[(4.5)]{BC} to be a sum of these $\mathscr{S}_{b,\xi}$ and $\mathscr{S}^*_{b,\xi}$, we can substitute these estimates into \cite[(4.33)]{BC} to yield
$$\mathscr{O}_{b}\ll \|\beta\|^2\|\nu\|^2 \left(b+ \frac{|\vartheta|A}{NM}\right)^\frac12ALN^{\frac34}X^{-\frac \eta 2}M^\varepsilon\bigg(  \frac{b^\frac14 N^\frac12 L^\frac12}{M^\frac12}+\frac{L^\frac52N}M+\frac{N^{\frac14}}{ A^\frac12X^{\frac \eta 2}}\bigg),$$
From \cite[(2.3)]{BC}, we have
$$\mathcal D_b=\mathscr{D}_b+\mathscr{O}_{b},$$
and from \cite[(2.1)]{BC}, we have
$$\mathcal{C}_b\ll ML^{-2+\varepsilon}\mathcal{D}_b,$$
and hence we can combine our estimates of $\mathscr{D}_b$ and $\mathscr{O}_{b}$ to change \cite[(5.1)]{BC} to
\begin{align*}\mathcal{C}_b &\ll \|\beta\|^2\|\nu\|^2 M^\varepsilon \left(1+ \frac{|\vartheta|A}{bNM}\right)^\frac12 \bigg(\frac{AM(bN)^{\frac 12}}{L^\frac12}+\frac{AM^2}{bLNX^{\eta}}+\frac{M^2}{LX^{\eta}}\\
&\hspace{10em}+ \frac{b^\frac34A M^\frac12 N^{\frac54}X^{-\frac \eta 2}}{ L^\frac12}+b^\frac12AL^\frac32N^{\frac74}X^{-\frac \eta 2}+\frac{b^\frac12 A^\frac12MN}{LX^{\eta}}\bigg).
\end{align*}
Choosing $L$ such that
$$b^\frac12AL^\frac32N^{\frac74}X^{-\frac \eta 2}\approx\frac{AM^2}{bLNX^{\eta}}+\frac{M^2}{LX^{\eta}}+\frac{b^\frac12 A^\frac12MN}{LX^{\eta}}, $$
we have
$$L=\frac{M^{\frac45}}{b^\frac35 N^{\frac{11}{10}}X^{\frac \eta 5}}+\frac{M^\frac45}{b^\frac15 A^\frac25N^{\frac7{10}}X^{\frac \eta 5}}+\frac{ M^\frac25}{A^\frac15N^{\frac3{10}}X^{\frac \eta 5}}+M^\varepsilon.$$
Plugging this into the estimate for $\mathcal C_b$ above then gives
\begin{align*}
\mathcal{C}_b &\ll \|\beta\|^2\|\nu\|^2M^\varepsilon \left(1+ \frac{|\vartheta|A}{bNM}\right)^\frac12 \bigg(AM(bN)^{\frac 12}+ b^\frac34A M^\frac12 N^{\frac54}X^{-\frac \eta 2}\notag\\
&\qquad+b^\frac12AN^{\frac74}X^{-\frac \eta 2}\bigg(\frac{M^{\frac45}}{b^\frac35 N^{\frac{11}{10}}X^{\frac \eta 5}}+\frac{M^\frac45}{b^\frac15 A^\frac25N^{\frac7{10}}X^{\frac \eta 5}}+\frac{ M^\frac25}{A^\frac15N^{\frac3{10}}X^{\frac \eta 5}}+M^\varepsilon\bigg)^\frac32\bigg)\notag\\
&\ll \|\beta\|^2\|\nu\|^2 M^\varepsilon \left(1+ \frac{|\vartheta|A}{bNM}\right)^\frac12 \bigg(AM(bN)^{\frac 12}+ b^\frac34A M^\frac12 N^{\frac54}X^{-\frac \eta 2}\notag\\
&\qquad+\frac{AM^{\frac65}N^{\frac1{10}}}{b^\frac2{5} X^{\frac{4\eta}{5}}}+b^\frac1{5}A^\frac25 M^\frac65N^{\frac7{10}}X^{-\frac{4\eta}{5}}+A^\frac7{10} b^\frac12M^\frac35N^{\frac{13}{10}}X^{-\frac{4\eta}{5}}+b^\frac12AN^{\frac74}X^{-\frac{\eta}{2}}\bigg).
\end{align*}
Following the methods of Section 6 of that paper to remove the square-free condition, one arrives at
\begin{align*}
 \mathcal{C}_1&\ll \|\beta\|^2\|\nu\|^2M^\varepsilon \left(1+ \frac{|\vartheta |A}{NM}\right)^\frac12 \\
&\times\hspace{0em}\big(AMN^{\frac 34}+ A M^\frac12 N^{\frac54}X^{-\frac \eta 2}+AM^{\frac65}N^{\frac1{10}-\frac{4\eta}{5}}+A^\frac25 M^\frac65N^{\frac7{10}}X^{-\frac{4\eta}{5}}+A^\frac7{10} M^\frac35N^{\frac{13}{10}}X^{-\frac{4\eta}{5}}+AN^{\frac74}X^{-\frac \eta 2}\big).
\end{align*}
As in the original paper, $AM^{\frac65}N^{\frac1{10}}X^{-\frac{4\eta}{5}}\ll AMN^{\frac 34}$ when $M\ll N^2$, while if $M\gg N^2$ then \cite[Theorem 5]{DFI} already gives a stronger bound of
$$\mathcal C_1\ll \|\beta\|^2\|\nu\|^2AM^{1+\varepsilon}.$$
Moreover, $A M^\frac12 N^{\frac54}X^{-\frac \eta 2}\ll AMN^{\frac 34}+AN^{\frac74}X^{-\frac \eta 2}$.  So
\begin{align*}
 \mathcal{C}_1&\ll \|\beta\|^2\|\nu\|^2M^\varepsilon \left(1+ \frac{|\vartheta |A}{NM}\right)^\frac12 \\
&\times\hspace{0em}\big(AMN^{\frac 34}+A^\frac25 M^\frac65N^{\frac7{10}}X^{-\frac{4\eta}{5}}+A^\frac7{10} M^\frac35N^{\frac{13}{10}}X^{-\frac{4\eta}{5}}+AN^{\frac74}X^{-\frac \eta 2}\big),
\end{align*}
and thus
\begin{align*}
 \mathcal{B}(M,N,A)&\ll \|\alpha\|\|\beta\|\|\nu\|M^\varepsilon \left(1+ \frac{|\vartheta |A}{NM}\right)^\frac12 \\
&\times\hspace{0em}\big(A^\frac 12M^\frac 12N^{\frac 38}+A^\frac15 M^\frac35N^{\frac7{20}}X^{-\frac{2\eta}{ 5}}+A^\frac7{20} M^\frac3{10}N^{\frac{13}{20}}X^{-\frac{2\eta}{5}}+AN^{\frac78}X^{-\frac \eta 4}\big),
\end{align*}
Again as in the original paper, we note that if $M\geq N$ then the parenthesized expression is dominated by $A^\frac 12M^\frac 12N^{\frac 38}+A^\frac7{20} M^\frac35N^{\frac7{20}-\frac{2\eta}{5}}$, whereas if $M\leq N$ then we can use B\'ezout's reciprocity theorem to switch $m$ and $n$, giving the same bounds with the $M$'s and $N$'s reversed.  So
\begin{align*}
 \mathcal{B}(M,N,A)&\ll \|\alpha\|\|\beta\|\|\nu\|M^\varepsilon \left(1+ \frac{|\vartheta |A}{NM}\right)^\frac12 \\
&\times\hspace{0em}\left(A^\frac 12\left(M^\frac 12N^{\frac 38}+M^\frac 38N^{\frac 12}\right)+A^\frac7{20} \left( M^\frac35N^{\frac7{20}}X^{-\frac{2\eta}{5}}+M^{\frac7{20}}N^{\frac35}X^{-\frac{2\eta}{5}}\right)\right).
\end{align*}
This completes the proof of Theorem \ref{BCsubdy}.

\section{Proof of Theorem \ref{FRch}: dispersion with subdyadic intervals}

Now, we turn to the Fouvry-Radziwiłł result.  Define
$$\Delta(\alpha, \beta, M, N, Q, a)=\sum_{\substack{q\sim Q \\ (a,q)=1}}\left|\mathop{\sum\sum}_{\substack{m \sim M \\ n \sim N \\ mn \equiv a \,(\mathrm{mod}\, q)}}\alpha_m \beta_n-\frac{1}{\varphi(q)}\mathop{\sum\sum}_{\substack{m \sim M \\ n \sim N \\ (mn,q)=1}}
\alpha_m \beta_n\right|.$$
Define $c_q=\pm 1$ or 0 to be such that if $(a,q)=1$ then $c_q=\pm 1$ with the sign matching the sign of the term inside of the absolute value above, while if $(a,q)>1$ then $c_q=0$.  We take a well-chosen smooth, compactly supported function $\psi$ with bounded derivatives (see \cite[Theorem 2.1]{FR2}), and we note that
\begin{gather}\label{Fourier1}\sum_{m\equiv a\pmod q}\psi\left(\frac mM\right)=\hat{\psi}(0)\frac Mq+\frac Mq\sum_{0<|h|\leq H}e\left(\frac{ah}{q}\right)\hat{\psi}\left(\frac{h}{q/M}\right)+O\left(M^{-1}\right),
\end{gather}
and
\begin{gather}\label{Fourier2}\sum_{(m,q)=1}\psi\left(\frac mM\right)=\frac{\varphi(q)}{q}\hat{\psi}(0)M+O\left(\tau(q)(\log 2M)^2\right).
\end{gather}
Letting $\mathcal L =\log 2X$, the authors use Cauchy-Schwarz to bound $\Delta^2$ as
$$  \Delta^2( \boldsymbol \alpha, \boldsymbol \beta, M, N, Q,a)  \ll MQ \mathcal L^{k^2-1}\, \Bigl\{ W(Q) -2 V(Q) +U(Q)\Bigr\},$$
where
\begin{align*}
  U(Q)&= \sum_{  (q,a)=1} \frac{\psi (q/Q)}{\varphi^2 (q)}\,\Bigl( \sum_{\substack{n \sim N \\ (n,q)=1}} \beta_n \Bigr)^2 \sum_{(m,q)=1} \psi \Bigl( \frac{m}{M} \Bigr),\\
  V(Q)& = \sum_{ (q,a)=1} \frac{\psi (q/Q)}{\varphi (q)}\,
  \Bigl( \sum_{\substack{n_1 \sim N \\ (n_1,q)=1}} \beta_{n_1} \Bigr)
  \Bigl( \sum_{\substack{n_2 \sim N \\ (n_2,q)=1}} \beta_{n_2} \Bigr) \sum_{  m\equiv a \overline{n_1}\bmod q } \psi \Bigl( \frac{m}{M} \Bigr),  \\
  W(Q) & =  \sum_{(q,a)=1} \psi (q/Q) 
  \Bigl( \sum_{\substack{n_1 \sim N\\ (n_1,q)=1}} \beta_{n_1} \Bigr)
  \Bigl( \sum_{\substack{n_2 \sim N\\ (n_2,q)=1}} \beta_{n_2} \Bigr) \sum_{\substack{ m\equiv a \overline{n_1}\bmod q \\ m\equiv a \overline{n_2}\bmod q}} \psi \Bigl( \frac{m}{M} \Bigr).
\end{align*}
We use the same estimates for $U(Q)$ and $V(Q)$ as were found in the original paper.  

For $W(Q)$, Fouvry and Radziwiłł apply (\ref{Fourier1}) in \cite[(26) and (27)]{FR} to find that 
$$W(Q)=\widetilde{W}^{MT}(Q)+\widetilde{W}^{Err1}(Q)+\widetilde{W}^{Err2}(Q)+O\left(MN^2Q^{-1}X^{\kappa-\frac \varepsilon 2}+X^{1+\kappa+\frac\varepsilon 2}\right),$$
with each of the three terms corresponding to the three terms in (\ref{Fourier1}) and where $\eta>0$ is arbitrarily small.  We will use the same estimates for $\widetilde{W}^{MT}(Q)$ and $\widetilde{W}^{Err2}(Q)$ as in the original paper and will focus our efforts on improving $\widetilde{W}^{Err1}(Q)$.

Let
$$H = M^{-1} Q X^\varepsilon.$$
We have
\begin{equation*}
  \widetilde{ W}^{\rm Err1} (Q) =  M\sum_{q} \frac {\psi (q/Q)}{q} \underset{\substack{n_1,  n_2 \sim N\\ n_1 \equiv n_2 \bmod q}}{\sum \sum} \beta_{n_1} \beta_{n_2}
   \sum_{0 < \vert h \vert \leq H} \hat \psi \Bigl( \frac{h}{q/M}\Bigr)\,e \Bigl( \frac{ah \overline{n_1}}{q}\Bigr):
 \end{equation*}
 Define
  \begin{equation}\label{decompn1n2}
  \begin{cases} (n_1, n_2 ) =d,\\
  n_1= d \nu_1, \ n_2 = d \nu_2, \  (\nu_1, \nu_2)=1, \\
  \nu_1= d_1\nu'_1\text { with } d_1 \mid d^\infty\text{ and } (\nu'_1, d)=1,
  \end{cases}
  \end{equation} 
and write
$$\nu_1-\nu_2=d\nu_1'-\nu_2=qr.$$
Then $1\leq r\leq R/d$, where $R$ is defined as in the original paper as
$$R=2NQ^{-1}.$$
Then
$$   \widetilde{W}^{\rm Err1} (Q) \ll X^{2\varepsilon } MQ^{-1}  \bigl\vert \, \mathcal W \,  \bigr\vert,$$
where
   \begin{equation*}
 \mathcal W = \sum_{1\leq \vert r \vert \leq R/d} \underset{\substack{dd_1\nu'_1, d\nu_2 \sim N \\ d_1\nu'_1\equiv  \nu_2 \bmod r}}
  {\sum \sum}\beta_{dd_1 \nu'_1} \beta_{d\nu_2}
   \frac{ \psi \bigl( (d_1\nu'_1-\nu_2)/ (rQ)\bigr)}{ (d_1\nu'_1-\nu_2)/ (rQ)}\\
    \sum_{0 <\vert h \vert \leq H}  \hat \psi \Bigl( \frac{h}{(d_1\nu'_1-\nu_2)/(rM)}\Bigr)
  e(\cdot ),
 \end{equation*}
   and $e(\cdot)$ is 
  $$ e(\cdot) = e \Bigl( \frac{ah\, \overline{dd_1 \nu'_1}}{(d_1 \nu'_1- \nu_2)/r} \Bigr),$$
  and  where the variables satisfy the following divisibility conditions: 
    $$  (d_1\nu'_1, \nu_2) =1,\ (\nu'_1, d) =1 \mbox{ and } (dd_1 \nu'_1 r, d_1\nu'_1 - \nu_2) =r.$$
For a given $r$, we now split the sum over $[N,2N]$ into subdyadic intervals of size $20Qr$.  For ease of notation, write
$$\mathcal N(U,r)=[N+20QrU,N+20Qr(U+1)).$$
Using Bézout's reciprocity formula as the authors do in the prelude to \cite[(39)]{FR2}, we then have
\begin{multline}\label{W1}
\vert  \mathcal W \bigr\vert \leq X^{6\varepsilon} \sum_{\substack{1\leq \vert r\vert \leq R/d \\ r \equiv a_1 \bmod{d d_1}}} \sum_{U_1=0}^{N/20Qr-1}\sum_{U_2=0}^{N/20Qr-1}\frac{1}{\phi(r)}\\
\Bigl\vert \sum_{\chi\pmod d}
\underset{\substack{dd_1\nu'_1\in \mathcal N(U_1,r) \\ d\nu_2\in \mathcal N(U_2,r)\\ \nu_1' \equiv a_2 \bmod d d_1 \\ \nu_2 \equiv a_3 \bmod d d_1}}
  {\sum \sum}\chi(d_1\nu'_1)\overline{\chi(\nu_2)}\beta_{dd_1 \nu'_1} \beta_{d\nu_2} \sum_{\substack{1\leq \vert h \vert \leq H \\ h \equiv a_4 \bmod d d_1}}
\Psi_r (h, \nu'_1, \nu_2)
 e \Bigl(  \frac{ah r \overline{dd_1\nu_2} }{\nu'_1} \Bigr)
\Bigr\vert,
\end{multline}
 where $\Psi_r$ is the differentiable function 
$$
\Psi_r (h, \nu'_1, \nu_2)=  \frac{ \psi \bigl( (d_1\nu'_1-\nu_2)/ (rQ)\bigr)}{ (d_1\nu'_1-\nu_2)/ (rQ)} \,  \hat \psi \Bigl( \frac{h}{(d_1\nu'_1-\nu_2)/(rM)}\Bigr)\, e \Bigl(\frac{ahr}{dd_1 \nu'_1 
( d_1\nu'_1 - \nu_2) } \Bigr),
$$
Note that if $|d_1\nu'_1-\nu_2|\geq 10rQd$ then $\Psi_r=0$.  So we can limit the sums over $U_1$ and $U_2$ to $|U_1-U_2|\leq X^\varepsilon$.


In order to remove the $\psi$ and $\hat \psi$ terms, we use partial summation over $n_1$, $n_2$, and $h$.  As in the original paper, we use Abel summation over the variables $\nu'_1$, $\nu_2$ and $h$, using the fact that for $0\leq \varepsilon_0,\, \varepsilon_1, \varepsilon_2 \leq 1$, we have the inequality   
\begin{equation}\label{bypart}
\frac{\partial^{\varepsilon_0+ \varepsilon_1+ \varepsilon_2}}{\partial h^{\varepsilon_0}\partial {\nu'_1}^{\varepsilon_1} \partial\nu_2^{\varepsilon_2}} \Psi_r (h, \nu'_1, \nu_2) \ll X^{50\varepsilon} \vert h\vert^{-\varepsilon_0} \, {\nu'_1}^{-\varepsilon_1} \, {\nu_2}^{-\varepsilon_2} \bigl(N/(rQ)\bigr)^{ \varepsilon_1 + \varepsilon_2 }.
\end{equation}
So
\begin{multline}\label{W2}
\vert  \mathcal W \bigr\vert \leq X^{\varepsilon} N^2Q^{-2}\sum_{\substack{1\leq \vert r\vert \leq R/d \\ r \equiv a_1 \bmod{d d_1}}} \frac 1{r^2\phi(r)} \sum_{U_1=0}^{N/20Qr-1}\sum_{\max\{0,U_1-X^\varepsilon\}\leq U_2\leq \min\{N/20Qr-1,U_1+X^\varepsilon\}}\left(\frac{Qr}{N}\right)^2\\
\Bigl\vert \sum_{\chi\pmod r}
\underset{\substack{dd_1\nu'_1\in \mathcal N(U_1,r) \\ d\nu_2\in \mathcal N(U_2,r)\\ \nu_1' \equiv a_2 \bmod d d_1 \\ \nu_2 \equiv a_3 \bmod d d_1}}
  {\sum \sum}\chi(d_1\nu'_1)\overline{\chi(\nu_2)}\beta_{dd_1 \nu'_1} \beta_{d\nu_2} \sum_{\substack{1\leq \vert h \vert \leq H \\ h \equiv a_4 \bmod d d_1}}
 e \Bigl(  \frac{ah r \overline{dd_1\nu_2} }{\nu'_1} \Bigr)
\Bigr\vert,
\end{multline}
We can combine some of these terms and restate the characters as congruence conditions, giving
\begin{multline}\label{W3}
\vert  \mathcal W \bigr\vert \leq X^{\varepsilon} \sum_{\substack{1\leq \vert r\vert \leq R/d \\ r \equiv a_1 \bmod{d d_1}}} \sum_{\substack{0< j<r \\ (j,r)=1}} \sum_{U_1=0}^{N/20Qr-1}\sum_{\max\{0,U_1-X^\varepsilon\}\leq U_2\leq \min\{N/20Qr-1,U_1+X^\varepsilon\}}\\
\Bigl\vert 
\underset{\substack{dd_1\nu'_1\in \mathcal N(U_1,r) \\ d\nu_2\in \mathcal N(U_2,r)\\ d\nu_1'\equiv \nu_2\equiv j\pmod r \\ \nu_1' \equiv a_2 \bmod d d_1 \\ \nu_2 \equiv a_3 \bmod d d_1}}
  {\sum \sum}\beta_{dd_1 \nu'_1} \beta_{d\nu_2} \sum_{\substack{1\leq \vert h \vert \leq H \\ h \equiv a_4 \bmod d d_1}}
 e \Bigl(  \frac{ah r \overline{dd_1\nu_2} }{\nu'_1} \Bigr)
\Bigr\vert,
\end{multline}
We now apply Theorem \ref{BCsubdy} to the inside expression.  Note that $n_1$ and $n_2$ are summed over intervals of size $\ll Qr$ but with only $1/r$ of the terms.  Hence, we can say that the intervals are of size $\ll Q$, which we can write as $N\cdot \left(\frac{Q}{N}\right)$.  Plugging in $\frac QN$ for $X^{-\eta}$ in the theorem, we have
\begin{align*}
\vert  \mathcal W \bigr\vert \leq &X^{\varepsilon} \sum_{\substack{1\leq \vert r\vert \leq R/d \\ r \equiv a_1 \bmod{d d_1}}} \sum_{\substack{0< j<r \\ (j,r)=1}} \sum_{U_1=0}^{N/20Qr-1}\sum_{U_1-X^\varepsilon\leq U_2\leq U_1+X^\varepsilon}\left(H^\frac 12Q\left(H^\frac 12N^\frac 78+H^\frac7{20} N^\frac{19}{20}\left(\frac{Q}{N}\right)^{\frac{2}{5}}\right)\right)\\
\ll &X^{\varepsilon} \sum_{\substack{1\leq \vert r\vert \leq R/d \\ r \equiv a_1 \bmod{d d_1}}} \frac{N\phi(r) }{Qr} \left(Q^2M^{-1}N^\frac 78+Q^{\frac{39}{20}}M^{-\frac{17}{20}} N^\frac{17}{20}\right)\\
\ll &X^{\varepsilon} \left(M^{-1}N^{\frac{23}{8}}+Q^{-\frac{1}{20}}M^{-\frac{17}{20}} N^{\frac{57}{20}}\right).\end{align*}
Hence,
$$   \widetilde{W}^{\rm Err1} (Q) \ll X^{\varepsilon } \left(Q^{-1}N^{\frac{23}{8}}+Q^{-\frac{21}{20}}M^{\frac{3}{20}} N^{\frac{57}{20}}\right)$$
Plugging this into the second equation in \cite[Section 4.5]{FR2} gives
\begin{multline*}
 \Delta^2 \ll MQ\mathcal L^{k^2-1} X^\varepsilon\Bigl\{  MN^2 Q^{-1} (\log 2N)^{-A} +  \\+   MN^2 Q^{-1} 
  + \bigl( M^\frac {3}{20} N^{\frac{57}{20}} Q^{-\frac{21}{20} } + N^{\frac{23}{8 }} Q^{-1}
  \bigr)
\Bigr\}, 
\end{multline*}
This changes the relevant equation in \cite[(43)]{FR2} to 
$$MQ \cdot M^\frac{3}{20} N^{\frac{57}{20}} Q^{-\frac{21}{20} }X^{\varepsilon}\ll M^2N^2,$$
which holds when
 \begin{gather}\label{14}
     Q> N^{34} X^{-17+ \varepsilon}.
\end{gather}
If we take $Q=X^{\frac 12+\varepsilon}$, we have
$$ X^{\frac{35}{68}-\varepsilon}> N.$$
This completes the proof of Theorem \ref{FRch}.

\end{document}